\documentclass[aps,prl,reprint,groupedaddress]{revtex4-1}
\usepackage{graphicx}
\usepackage{natbib}

\begin{document}


\title{Structural Susceptibility and Separation of Time Scales in the van der Pol Oscillator}
\author{Ricky Chachra}
\author{Mark K. Transtrum}
\author{James P. Sethna}
 \email[]{sethna@lassp.cornell.edu}

\affiliation{Department of Physics, Laboratory of Atomic and Solid State Physics\\Cornell University, Ithaca, NY 14853 USA}




\date{\today}

\begin{abstract}

We use an extension of the van der Pol oscillator as an example of a system with multiple time
scales to study the \textit{susceptibility} of its trajectory to polynomial perturbations in the dynamics.
A striking feature of many nonlinear, multi-parameter models is an apparently inherent 
insensitivity to large magnitude variations in certain linear combinations of parameters. 
This phenomenon of ``sloppiness'' is quantified by calculating the eigenvalues of the 
Hessian matrix of the least-squares cost function. These   
typically span many orders of magnitude. The van der Pol system is no exception: 
Perturbations in its dynamics show that most directions in parameter space weakly 
affect the limit cycle, whereas only a few directions are stiff. 
With this study we show that separating the time scales in the van der Pol system
leads to a further separation of eigenvalues. Parameter combinations which 
perturb the slow manifold are stiffer and those which solely affect the jumps in the 
dynamics are sloppier.

\end{abstract}

\pacs{}

\maketitle


\section{Introduction}

In this manuscript, we analyze the sensitivity of a multiple time scales dynamical system to perturbative 
changes in its evolution laws. Rather than utilizing the traditional means of examining the \textit{structural 
stability} for probing qualitative changes to the attractor as a response to perturbations, 
we study the \textit{structural susceptibility} for quantifying the effects of the perturbations 
on the time series~\footnote[1]{We employ the word \textit{structural} in the same context as its usage in 
dynamical systems literature on \textit{structural stability}.  The word \textit{susceptibility} is inspired from physics 
wherein it is a measure of response to a perturbation (such as an applied external field) quantified by the second-derivative 
of the free energy w.r.t. parameters. Since cost is analogous to free energy (in that both are minimized), 
it is natural to call the response to perturbations in dynamics, also quantified via second derivatives, 
as \textit{structural susceptibility}}. 
More specifically, we ask how sensitive is the dynamical system 
$d\textbf{z}/dt=\textbf{f}(\textbf{z})$ to infinitesimal changes of the form 
$d\textbf{z}/dt=\textbf{f}(\textbf{z}) + \textbf{a} \cdot \textbf{g}(\textbf{z})$, for a family of perturbations $\textbf{g}(\textbf{z})$
when the parameters $\textbf{a}\rightarrow\textbf{0}$.	

This report introduces the new concept of ``structural susceptibility'' in dynamical systems, 
and is an outgrowth of our group's previous  work on ``sloppiness'' in multiparameter 
systems wherein we have found that 
data-fitting in a number of nonlinear, multiparameter models is only sensitive to a few directions 
in parameter space at the best fit~\cite{BrownSethna, Guntenkunst1, MarkLongPaper}. The key difference 
between studying 
sloppiness and structural susceptibilities is that in the former, the parameters are intrinsic to  
the system, i.e., there are no externally introduced changes in their evolution laws.
Nonetheless, the methodology we have developed for studying sloppy models is also suited for studying 
structural susceptibilities of dynamical systems. Our approach cleanly isolates and ranks the 
directions in parameter space in order of relevance to observed behavior, 
and has previously led us to suggest improvements in experimental design~\cite{Optimal}, 
extract falsifiable predictions from experiments~\cite{Falsify}, 
and develop faster minimization algorithms~\cite{MarkLM}. Others have developed our
ideas to suggest further improvements in experimental design~\cite{Apgar} and parameter 
estimation~\cite{Secrier}, to quantify robustness to parameter variations~\cite{Dayarian}, and to 
set confidence regions for predictions in multiscale models~\cite{Multiscale}. In this paper, we bring similar 
ideas together to analyze sensitivities of time series to perturbations in dynamical systems. 
 
We demonstrate the utility of our approach with application to a dynamical system with two time scales---
the van der Pol oscillator~\cite{VDP} which is a single parameter system and hence not amenable to 
sloppy model analysis. Instead, by choosing appropriate perturbations $\textbf{g}(\textbf{z})$, we  
calculate the susceptibility of its dynamics: We make perturbations on the attractor, and then systematically 
increase the separation of time scales in its dynamics to show how it can generally enhance the sloppiness in 
nonlinear systems.

\section{Multiple Time Scale Dynamics}

Multiple time scales are often found in the solutions of dynamical systems~\cite{Jones}.
Broadly speaking, the defining criterion of these models is that the trajectory 
of one or more phase variables has 
an identifiable fast piece such as a jump or a pulse and a slow piece where the value of the 
variable doesn't change quickly~\cite{GrasmanBook}.  In two dimensions, these systems are 
commonly studied in the contexts of slow-fast vector fields written as: 
\begin{equation}
	\begin{array}{rl}
	\epsilon \dot{x}&=X(x, y, \epsilon), \\
	\dot{y}&=Y(x,y, \epsilon)
	\end{array} 
	\label{slowfast}
\end{equation}
where the parameter $\epsilon>0$ is small and dot indicates derivative with respect to time $t$. 
For ${\cal O}(1)$ functions $X$ and $Y$, and $X\neq0$: $\dot{x}={\cal O}(1/\epsilon)$ 
and $\dot{y}={\cal O}(1)$, so that $\epsilon$ is the ratio of time scales in the system. 
On one extreme, the singular limit $\epsilon = 0$ corresponds to a differential 
algebraic system $X=0,\ Y=\dot{y}$ where the solutions of $X = 0$ comprise the 
``critical manifold'' close to which the flow in phase space is slow (the ``slow manifold''). Similarly, $\epsilon=1$ 
corresponds to a limit where there is no separation of time scales, with a crossover at intermediate values
of $\epsilon$.

Originally introduced in 1927, the van der Pol equation, $\ddot{x} - \mu(1-x^2)\dot{x} +x=0$, 
is a well-studied example of a second-order, nonlinear system with multiple time scales in its solution. 
Using the Li\'{e}nard 
transformation $y=x-x^3/3-\dot{x}/\mu$, and redefining time $t\rightarrow t\mu$, the equation 
can be written as a two dimensional system~\cite{GrasmanBook, Strogatz} given by:
\begin{equation}
\begin{array}{rl}
	\mu^{-2} \dot{x}&=x - \frac{x^{3}}{3} - y, \\
	\dot{y}&=x, 
\end{array}
\label{vdpeqns}
\end{equation}
which has the same form as (\ref{slowfast}) with $\epsilon=\mu^{-2}$. The global attractor of this dynamical system is 
a structurally stable limit cycle with two time scales~\footnote[2]{Incidentally, the set of equations (\ref{vdpeqns}) 
can also be considered a special case of the FitzHugh-Nagumo model 
(\textit{See} R. FitzHugh. Impulses and Physiological States in theoretical models of nerve propagation \textit{Biophys J.}, 1(6):445, 1961) introduced three decades later 
as simplification of the Hodgkin-Huxley equations of neuronal spikes in the squid giant axons, and 
is sometimes referred to as the Bonhoeffer-van der Pol model}. 

\begin{figure}[h]
 \begin{center}
  \includegraphics[scale=0.8, trim= 115mm 160mm 13mm 14mm, clip, angle=90]{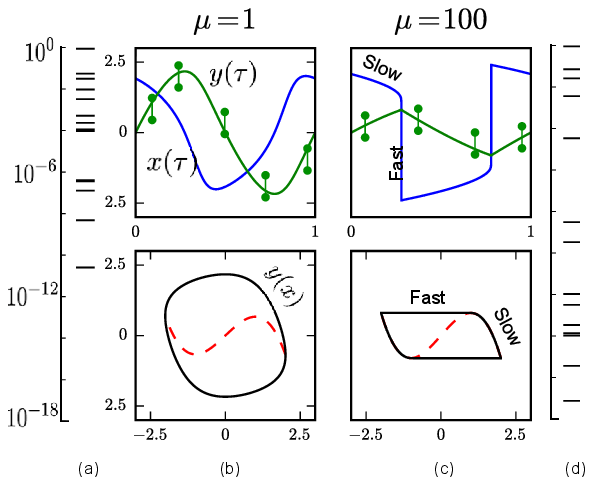}
  \caption{\label{fig:trajslow} (a) Eigenvalues of the Hessian matrix of the cost of fitting at $\mu=1$, and 
  (b, c top row) one period of time series x($\tau$) and y($\tau$), for $0<\tau<1$, 
  are shown for $\mu=1$ and $\mu=100$ as function of time along with schematic error bars for 
  the data-fitting of the trajectory of variable $y$. (d) Eigenvalues for $\mu=100$ are shown on the right. 
  As $\mu\rightarrow\infty$, the orbit collapses onto the critical manifold with the trajectory spending most 
  of its time on the slow manifold and vanishingly short on the jumps. Also shown in (b, c bottom row) 
  is the orbit in $xy$ plane (solid line) and the critical manifold (dashed line). }
   \end{center}	
\end{figure}

The van der Pol system provides a convenient way 
to separate time scales by varying $\mu$:
Small values of $\mu$ in the van der Pol system correspond to a small separation of 
time scales. It can be shown that the trajectory approaches that of the 
harmonic oscillator as $\mu\rightarrow0$~\cite{VDP}.
At large values of $\mu$, the system shows a separation of   
time scales which increases with increasing $\mu$. As shown in figure~\ref{fig:trajslow}
(b, c), with increasing $\mu$, the trajectory of $x$ separates into a slow part that lies ${\cal O}(exp(\mu^{-2}))$ 
close to the phase space curve given by $\dot{x} = 0$, \textit{i.e.} the critical manifold $y = x 
- x^3/3$, and a fast part which connects the two branches of the slow flow.
Likewise, the separation of time scales in $y$ are associated with the increasing sharpness
of the kink in its trajectory. 

The fact that with an increasing separation of time scales the trajectory 
spends an increasing amount time on the slow manifold and a decreasing amount of time on the 
jumps has important implications for fitting parameters to time series data of the van der Pol system. 
With increasing scale separation, one expects that the cost of fitting 
will be decreasingly sensitive to changes in the jumps of the trajectory as they get progressively 
shorter in duration. 

\section{Sloppiness in Nonlinear Fits}

In this section, we discuss the concepts of sloppiness and structural susceptibility in more detail with examples as a 
prelude to our calculations. For time series $z(t, \textbf{a})$, a least-squares fit to data $d_i$ minimizes a cost $C=\frac{1}{2}\sum_i (z(t_i, \textbf{a})-d_i)^2/\sigma_i^2$ in the space of system parameters which are collectively denoted as \textbf{a}.
Our discovery of sloppiness is essentially that the eigenvalues of the Hessian of the cost with respect to parameters,  
${\cal H}_{\alpha\beta} = \partial^2C/\partial a_\alpha \partial a_\beta$, at the best fit span many orders of magnitude.
The larger and smaller eigenvalues correspond to stiffer and sloppier directions respectively. 
For concreteness, consider a time series of a multi parameter model, such as the one denoted by $y(\tau)$ in figure~\ref{fig:trajslow}(b, top row). The error bars schematically show the least-squares fit of $y(\tau)$ and the sidebar (figure~\ref{fig:trajslow}(a)) shows the eigenvalues of the corresponding Hessian matrix. Note the broad range of eigenvalues ($\sim 10^{11}$, corresponding to a factor of almost a million in parameter range)--- a feature that turns out to be typical in nonlinear fits. 

Another vivid example of sloppiness in nonlinear models is provided by the well-established formalism behind the 
characterization of the sensitivities of initial conditions using Lyapunov exponents~\cite{Katok}.  
Consider $d\textbf{z}/dt=\textbf{f}(\textbf{z})$ as a model whose parameters are the initial conditions 
$a_\alpha=z_\alpha(0)$ and whose predictions are the final positions $z_i(t)$ at time $t$. At the best fit, 
${\cal H}_{\alpha\beta}=  (J^TJ)_{\alpha\beta}$ where $J_{i\alpha}=\partial z_i(t)/\partial z_\alpha(0)$
is the Jacobian of the sensitivities to perturbations in the initial conditions. 
The Lyapunov exponents, which are defined to be the eigenvalues $\ell_n$ of $\textbf{L} = \lim_{t\to\infty} 1/(2t) \log(J^TJ)$,
utilize the same Hessian we would use in calculating the sloppy model eigenvalues $\lambda_n=exp(2t\ell_n)$. 
The roughly equal spacing of Lyapunov exponents naturally explains 
both the exponentially broad range of sloppy model exponents and the roughly equal spacing 
of $\log(\lambda_n)$ for a model with initial conditions as parameters.

Instead of the sensitivities with respect to the initial conditions or other intrinsic parameters, we focus here on the 
sensitivity of the dynamics to changes in the dynamical evolution laws. Therefore, for the remainder of this paper we will be interested in a cost function that measures the square of the distance between two time series for the system $d\textbf{z}/dt=\textbf{f(z)}+\textbf{a}\cdot \textbf{g}(z)$--- one with perturbations--- $\textbf{z}(t, \textbf{a}\rightarrow\textbf{0})$, and the other one without, \textit{i.e.}, $\textbf{z}(t, \textbf{a=0})$
\begin{eqnarray}
C = \frac{1}{2} \int_0^T \! || \textbf{z}(t, \textbf{a} \rightarrow \textbf{0}) - \textbf{z}(t, \textbf{a}=\textbf{0})||^2 \, \mathrm{d}t 
\label{suscost}
\end{eqnarray}
with the perturbing terms $g_i(z)$ giving a power series in the components of $\textbf{z}$. 
Further in the manuscript, we will use this form of the cost to compute the susceptibility of the van der Pol system and show how sloppiness is enhanced by increasing separation of time scales in the van der Pol equations. This is in essence captured by 
figure~\ref{fig:trajslow}(a \& d) where we show that an increase in the van der Pol parameter $\mu$ from $1$ to $100$ 
produces roughly a million-fold increase in the spread of eigenvalues.

\section{Susceptibility of van der Pol system}


We perturb the van der Pol system in (\ref{vdpeqns}) by adding a series of additional 
terms. There is a long tradition in dynamical systems of studying
equations of motion of polynomial form~\cite{Katok, Gucken}; indeed, the theory of normal forms
suggests that general dynamical systems, even at bifurcations, can be
generically mapped into a polynomial form by a nonlinear but smooth change
of variables.  Adding extra polynomial terms is routinely done to `unfold' the qualitative behavior 
near bifurcations~\cite{Unfold}. Here we focus on quantitative changes far from bifurcations. 
In choosing our perturbations, we must cut off the polynomials at some order. There are two ways in which 
we specialize our general susceptibility analysis to the two time scale, 
periodic limit cycle of the van der Pol system. First, we choose the family of perturbations of 
order $3N$ as follows:
\begin{equation}
\begin{array}{rl}
	\mu^{-2} \dot{x}&=x - \frac{x^{3}}{3} - y + 
	\sum_{m+n\leq N} a_{m,n}(x - \frac{x^{3}}{3} - y)^{m}x^{n} \nonumber 
	\\
	\dot{y}&=x. \nonumber
 \end{array}
 \label{vdppert}
\end{equation}
This choice has two noteworthy features---
(a)~We have grouped the polynomial perturbations so that, for $m\ne0$ they
vanish on the critical manifold, $y=x-x^3/3$. That is, the parameters $a_{m,n}$ with $m\ne0$
do not significantly affect the dynamics on the slow manifold; we call these
the ``fast parameters'' and correspondingly the $a_{0,n}$ are
``slow parameters''.  The parameter $a_{1,0}$ duplicates $\mu$ to the
 same effect as changing the period, and we omit it. 
 Surely, the eigenvalue spectrum of the general polynomial
expansion, $a_{m,n}x^my^n$, behaves qualitatively similarly to the one presented here
but our parametrization greatly simplifies the analysis of the eigenvector perturbations. 
(b)~We only perturb the $\dot x$ equation. Our choice corresponds to a 
general expansion of a second-order equation, with the acceleration
$\ddot y = \dot x$ written as a polynomial in the position $y$ and
velocity $\dot y = x$. Perturbing both equations produces similar behavior. 

Second, we modified the cost to focus on the limit cycle of the van der Pol system in two ways--- (a) by rescaling 
all trajectories in our analysis so that they have the same unit period, and (b) by changing the initial condition so that it lies on the perturbed orbit and that the perturbed and unperturbed orbits are in phase with each other~\footnote[3]{Perturbations distort the dynamics so that the attractor and its period change. We addressed these issues by setting the periods to unity, and by moving the initial conditions to the new attractor to remove any transients. Alternatively, if we fit data over many periods without making the said changes, the parameter combinations determining the period and phase would become stiff modes in our dynamics.}.
When we correct the period $T$ by $\delta T$, initial conditions $\textbf{y}_0$ by
 $\delta \textbf{y}_0$, and do an overall rescaling of the time variable 
 $t \rightarrow \tau T$,  the cost functional for the time series of $y(\tau)$ at each $\mu$ 
 takes the following form:
\begin{eqnarray}
	C(\mu) = \frac{1}{2} \int_0^1 \! [ y(\tau, \textbf{a}+\delta \textbf{a},
	\textbf{y}_0 + \delta \textbf{y}_0, T+\delta T)-
	\nonumber \\ y(\tau, \textbf{a}, \textbf{y}_0, T) ]^2 \, \mathrm{d}\tau 
\end{eqnarray} 
In principle, changes in both time series, $x(\tau)$ and $y(\tau)$ could be 
incorporated in the cost function, but we get qualitatively similar 
results by keeping either or both variables. Choosing to measure changes only in $y(\tau)$  
corresponds again to studying the second-order equation for $\ddot{y}$ as an expansion in $y$ and $\dot{y}$.

\begin{figure}[h]
 \begin{center}
  \includegraphics[scale=1.0, trim= 150mm 175mm 13mm 14mm, clip, angle=90]{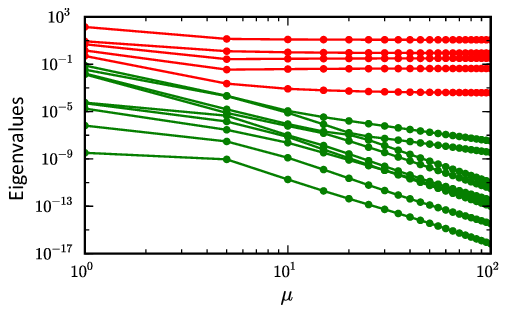}
  \caption{\label{fig:vdpevals} Eigenvalues of Hessian matrix are shown here as a 
  function of $\mu$. The range $1 \leq\mu\leq 100$ corresponds to a ratio of time
  scales $1\leq\epsilon\leq10000$. The five largest eigenvalues correspond to stiff 
  directions in the parameter space: these directions perturb the slow manifold. The
  remainder affect the transient part of the trajectory which becomes smaller with 
  an increasing separation of time scales and hence these directions are 
  decreasingly relevant.}
 \end{center}	
\end{figure}

The susceptibilities are still given by the Hessian matrix at the best fit ($\textbf{a=0}$): 
\begin{eqnarray}
	{\cal H}(\mu)_{\alpha\beta} = \frac{\partial^2 C(\mu)}{\partial a_\alpha \partial a_\beta} 
\end{eqnarray}
which can be written out more completely as:
\begin{eqnarray}
	{\cal H}(\mu)_{\alpha\beta} = \int_0^1 \! \left( \frac{\partial y}{\partial a_\alpha}
	+ \frac{\partial y}{\partial \textbf{y}_0}\frac{\partial \textbf{y}_0}{\partial a_\alpha}
	+ \frac{\partial y}{\partial T}\frac{\partial T}{\partial a_\alpha}\right) \nonumber 
	\\ \times \left( \frac{\partial y}{\partial a_\beta}
	+ \frac{\partial y}{\partial \textbf{y}_0}\frac{\partial \textbf{y}_0}{\partial a_\beta}
	+ \frac{\partial y}{\partial T}\frac{\partial T}{\partial a_\beta} \right)
   \, \mathrm{d}\tau  \nonumber
\end{eqnarray}
Here, each of the two terms in the integral is to be interpreted as a Jacobian matrix, 
a mapping from the finite dimensional parameter space to the infinite dimensional data space:
\begin{eqnarray}
J_{\tau\alpha} =  \frac{\partial y(\tau)}{\partial a_\alpha}
	+ \frac{\partial y(\tau)}{\partial \textbf{y}_0}\frac{\partial \textbf{y}_0}{\partial a_\alpha}
	+ \frac{\partial y(\tau)}{\partial T}\frac{\partial T}{\partial a_\alpha}   
\label{Jacobian}
\end{eqnarray}
The sensitivity trajectories in the Jacobian, $\partial y/ \partial a_\alpha $, $\partial y/ \partial \textbf{y}_0$,  
and $\partial y/ \partial T $, were computed using the open source SloppyCell package~\cite{SloppyCell1, SloppyCell2}. 
The expressions for the time invariant quantities,  $\partial \textbf{y}_0 / \partial a_\alpha $ and 
$\partial T/ \partial a_\alpha $, were 
found by enforcing periodicity of the perturbed time series denoted by $\textbf{y}(\tau) \equiv (x(\tau), y(\tau))$ 
as follows: 
\begin{eqnarray}
\textbf{y}(\tau=0, \textbf{a}+\delta \textbf{a}, \textbf{y}_0 + \delta \textbf{y}_0, T+\delta T) =  \nonumber \\
  \textbf{y} (\tau=1, \textbf{a}+\delta \textbf{a}, \textbf{y}_0 + \delta \textbf{y}_0, T+\delta T),  \nonumber
\end{eqnarray}
Taylor expansion of both sides of the previous equation leads to a vector equation:
\begin{eqnarray}
\delta \textbf{y}_0 = \frac{\partial \textbf{y}}{\partial T}\bigg|_{\tau=1} \delta T + \frac{\partial \textbf{y}}{\partial \textbf{a}}
\bigg|_{\tau=1} \delta \textbf{a}  
+    \frac{\partial \textbf{y}}{\partial \textbf{y}_0} \bigg|_{\tau=1} \delta \textbf{y}_0,  \nonumber 
\end{eqnarray}
from which both constants can be computed following the convention that the component denoting
the change in initial condition of $y(\tau)$ in $\delta \textbf{y}_0$ is set to zero. Now with the Jacobian calculated, 
the Hessian at best fit is simply ${\cal H} = J^TJ $.  

\begin{figure*}[ht]
 \begin{center}
  \includegraphics[scale=1.1, trim= 135mm 105mm 15mm 15mm, clip, angle=90]{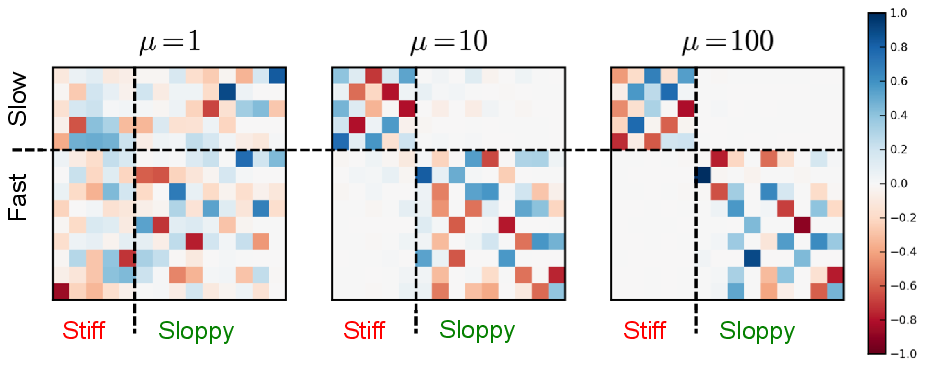}
  \caption{\label{fig:vdpeves} Hessian eigenvectors are shown for $\mu=1,\ 10,$ and $100$.
  Each colored small square shows the magnitude of an eigenvector component (the scale bar shown on the right). 
  Eigenvectors for each $\mu$ are sorted so that the stiffer ones appear on the left; individual components are sorted so that 
  ``slow parameters'' appear on the top. Note that with increasing $\mu$, the stiff and sloppy eigenvectors separate by parameters:     The stiff eigenvectors only have projections along the slow parameters which perturb the slow manifold, whereas
  the sloppy directions have projections along the fast parameters which mainly perturb the jumps.}
 \end{center}	
\end{figure*}

\subsection{Eigenvalues and Eigenvectors}

\begin{figure*}[ht]
 \begin{center}
  \includegraphics[scale=1.0, trim= 88mm 90mm 15mm 12mm, clip, angle=90]{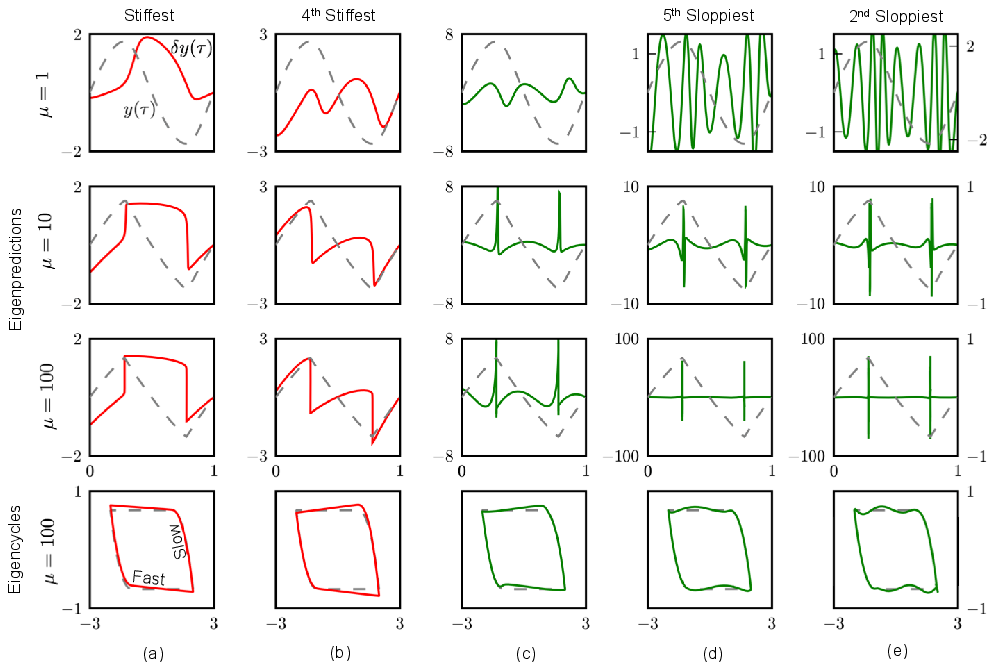}
  \caption{\label{fig:pertraj} Top three rows: Eigenpredictions $\delta y_k$ for $k=0,\ 3,\ 6,\ 9,\ 12$ at $\mu=1,\ 10\ \&\ 100$ are shown in solid red lines for stiff modes and solid green for sloppy modes. These curves show the response of perturbations if the parameters are changed infinitesimally along the Hessian eigenvectors: A parameter change of norm $\epsilon$ along eigendirection $n$ will change the trajectories by $\lambda_n \epsilon$ times the eigenpredicton. Dashed lines show unperturbed van der Pol solution for comparison ($y$ scale on the right hand side). As the time scales separate, the amplitudes of the sloppiest eigenpredictions increase (roughly in proportion to $\mu$) getting increasingly concentrated at the jumps. Bottom row shows the eigencycles for $\mu=100$ in solid red lines corresponding to the perturbations in row 3 (i.e. the new limit cycle for a perturbation of strength $\epsilon\sim 1/\lambda_n$. These curves show how the van der Pol orbit changes with perturbations along the Hessian eigenvectors. Both the stiff and the sloppy modes change the orbit at the jumps (occurring at the extrema in the dashed lines); the stiff modes also change behavior at the slow manifold, whereas the sloppy modes only affect the jumps.}
 \end{center}	
\end{figure*}

We computed the Hessian matrix given by the previous equation at multiple values of $\mu$. 
The spread of eigenvalues (figure~\ref{fig:vdpevals}) increases as a function of $\mu$ confirming that  
sloppiness increases with an increasing separation of time scales. Not surprisingly~\footnote[5]{We 
understand this as sloppiness as arising due to the generalized interpolation argument~\cite{MarkLongPaper}.}, 
for $N=4$, the $14$ eigenvalues for $\mu = 1$ already span $11$ orders of magnitude, while for 
$\mu = 100$, we observe that the stiffest eigenvalue is $18$ orders of magnitude 
larger than the smallest one--- the spread increases by $10^7$ when $\mu$ increases to $100$.

Taken together with the eigenvectors shown in 
figure~\ref{fig:vdpeves}, some interesting facts come to light:
Figure~\ref{fig:vdpevals} shows that with increasing $\mu$, the eigenvalues separate into two clusters of closely 
related decay exponents.  The largest $N$ eigenvalues approach  
constants. The other eigenvalues decay 
with power laws: two modes with exponents between $-2$ and $-3$ and the remaining 
with exponents between $-5$ and $-6$.  Similarly, figure~\ref{fig:vdpeves} shows that the eigenvectors also separate into 
two groups with increasing $\mu$: The stiffest directions are linear combinations 
of the slow parameters whereas the sloppy directions are comprised of other parameters
as expected.

We can understand the effect of perturbations in parameter combinations given 
by the eigenvectors (called \textit{eigenparameters}) $\hat{\textbf{e}}_k$ 
more clearly by observing their behavior 
in the data space. The Jacobian transformation of (\ref{Jacobian}) projects the eigenvectors
 to the data space: $\delta y_k = J \cdot \hat{\textbf{e}}_k / \sqrt{\lambda_k}$
where ${\lambda_k}$ corresponds to the $k^{th}$ largest eigenvalue. Defined this way, 
these data space vectors, called \textit{eigenpredictions}~\cite{MarkLongPaper}, $\delta y_k$, are also orthonormal. Alternatively,  
the eigenpredictions are the left singular vectors in the singular value decomposition of the Jacobian 
(i.e. they are the columns of the unitary matrix $U$ in $J = U\Sigma V^T$~\cite{NR}) 
As shown in figure~\ref{fig:pertraj} for $\mu=1,\ 10\ \&\ 100$ (top three rows), we learn from the
eigenpredictions that the stiff modes affect behavior both along the slow manifold and at the jumps. 
Moreover with increasing $\mu$, as the eigenvalues associated with the stiff directions 
approach constants (figure~\ref{fig:vdpevals}), so do the stiff eigenpredictions 
(figure~\ref{fig:pertraj} rows 2,\ 3 columns (a) and (b)). The sloppy modes 
on the other hand, affect the dynamics on the jumps only. The maximum amplitudes of the (normalized) sloppiest 
eigenpredictions appears to increase roughly proportional to $\mu$ (corresponding to the jump duration 
of $\sim\mu^{-2}$). In the limit, these become 
$\delta$-functions and derivatives concentrated at the jumps. 
Figure~\ref{fig:pertraj} (bottom row) also shows the limit cycles (\textit{eigencycles}) with 
eigenparameter perturbations as phase space trajectories $(x,  y + \eta\ \delta y_k)$ for 
small $\eta$.

\section{Discussion} 

In this paper, we have introduced a formalism we call ``structural susceptibility'' for analyzing the quantitative 
dependence of dynamical systems to perturbations of the equations of motion. It is a generalization
of `unfolding' methods of bifurcation theory and the Lyapunov exponents governing the dependence on initial 
conditions, and exposes the ubiquitous presence of broad range of sloppy eigendirections in parameter space--- largely
unimportant to the dynamics. We used this method to study the role of time scale separation in enhancing the sloppiness
of the susceptibility spectrum in the particular case of the van der Pol oscillator. 

By extending the framework of our sloppy model analysis to systems
where changes in evolution laws are to be studied, our method offers a simple way to 
calculate the effects of broad classes of perturbations. 
By studying the structural susceptibility of a dynamical system with two time scales, 
the analysis presented here showed that sloppiness of nonlinear systems
is enhanced by separation of time scales in the dynamics.  With increasing separation of 
time scales in the van der Pol oscillator, the trajectory spends an increasing 
amount of time on the slow manifold and a vanishingly small amount of time in the transition 
region. The cost of perturbations is integrated over time and therefore we 
are unsurprised that the perturbations that change the slow manifold will accrue the 
most cost and therefore manifest as stiff modes of the Hessian matrix. 
The remaining directions are sloppy as they only affect the behavior 
at the jumps or the fast pieces. These perturbations manifest as $\delta$-functions and their 
derivatives--- significantly affecting the phase-space trajectory, but over only the fast times 
asymptotically ignored in the least-squares cost. It remains a challenge to connect separation of time scales
to parameter sensitivity  in more complicated systems, but the analogy of the van der Pol system's behavior 
with other nonlinear physical systems of interest is clear. 

Many important dynamical systems have multiple time scales in their solutions: 
examples include models in neuroscience  
(such as Hodgkin-Huxley model), systems biology or chemical reaction systems 
(such as protein network models), and in engineering (such as models of combustion,
lasers, locomotion, etc.). Our analysis suggests that any system with multiple time scales should become sloppier 
as the scales separate for
the same reasons as we found in the van der Pol: Some parameter combinations 
will only affect the fast dynamics, and lead to sloppy modes. Perturbations which affect the slow dynamics will 
presumably accrue more cost and be stiff. 

More broadly, the sloppiness exposed by our structural susceptibility analysis has clear implications for attempting to 
reconstruct the equations of motion from experimental data~\cite{Lipson} because the true dynamics along any
sloppy eigendirection will be relatively poorly determined. This discovery has already influenced work on 
experimental design optimization: estimating parameters is challenging~\cite{Apgar, Comment}, but extracting predictions 
without constraining parameters is straightforward~\cite{Falsify}. We further anticipate that the concept of structural susceptibility will be useful for studying systems with chaos, bifurcations and phase transitions; quantifying the unfoldings of these systems 
may also be useful for gaining a deeper understanding of the phenomena they model.     

%
%
%

\section{Acknowledgments}

We thank Stefanos Papanikolaou for helpful conversations leading 
us to think about perturbing the slow manifold and about the analogues to thermodynamic susceptibilities, 
and John Guckenheimer for valuable insights and discussions regarding our calculations. 
Support from NSF grant DMR 1005479 is gratefully acknowledged.

\bibliographystyle{unsrt}
\bibliography{SloppyRefs}

\begin{thebibliography}{10}

\bibitem{Note1}
We employ the word \protect \textit {structural} in the same context as its
  usage in dynamical systems literature on \protect \textit {structural
  stability}. The word \protect \textit {susceptibility} is inspired from
  physics wherein it is a measure of response to a perturbation (such as an
  applied external field) quantified by the second-derivative of the free
  energy w.r.t. parameters. Since cost is analogous to free energy (in that
  both are minimized), it is natural to call the response to perturbations in
  dynamics, also quantified via second derivatives, as \protect \textit
  {structural susceptibility}.

\bibitem{BrownSethna}
Kevin~S. Brown and James~P. Sethna.
\newblock Statistical mechanical approaches to models with many poorly known
  parameters.
\newblock {\em Phys. Rev. E}, 68:021904, Aug 2003.

\bibitem{Guntenkunst1}
Ryan~N Gutenkunst, Joshua~J Waterfall, Fergal~P Casey, Kevin~S Brown,
  Christopher~R Myers, and James~P Sethna.
\newblock Universally sloppy parameter sensitivities in systems biology models.
\newblock {\em PLoS Comput Biol}, 3(10):e189, 10 2007.

\bibitem{MarkLongPaper}
Mark~K. Transtrum, Benjamin~B. Machta, and James~P. Sethna.
\newblock Geometry of nonlinear least squares with applications to sloppy
  models and optimization.
\newblock {\em Phys. Rev. E}, 83:036701, Mar 2011.

\bibitem{Optimal}
Fergel~P. Casey, D.~Baird, Q.~Feng, R.N. Gutenkunst, J.J. Waterfall, C.R.
  Myers, K.S. Brown, R.A. Cerione, and J.P. Sethna.
\newblock Optimal experimental design in an epidermal growth factor receptor
  signalling and down-regulation model.
\newblock {\em IET Systems Biology}, 1(3):190--202, 2007.

\bibitem{Falsify}
Ryan~N. Gutenkunst, Fergal~P. Casey, Joshua~J. Waterfall, Christopher~R. Myers,
  and James~P. Sethna.
\newblock Extracting falsifiable predictions from sloppy models.
\newblock {\em Annals of the New York Academy of Sciences}, 1115(1):203--211,
  2007.

\bibitem{MarkLM}
Mark~K. Transtrum and James~P. Sethna.
\newblock Improvements to the levenberg-marquardt algorithm for nonlinear
  least-squares minimization.

\bibitem{Apgar}
Joshua~F. Apgar, David~K. Witmer, Forest~M. White, and Bruce Tidor.
\newblock Sloppy models{,} parameter uncertainty{,} and the role of
  experimental design.
\newblock {\em Mol. BioSyst.}, 6:1890--1900, 2010.

\bibitem{Secrier}
Maria Secrier, Tina Toni, and Michael P.~H. Stumpf.
\newblock The {ABC} of reverse engineering biological signalling systems.
\newblock {\em Mol. BioSyst.}, 5:1925--1935, 2009.

\bibitem{Dayarian}
Adel Dayarian, Madalena Chaves, Eduardo~D. Sontag, and Anirvan~M. Sengupta.
\newblock Shape, size, and robustness: Feasible regions in the parameter space
  of biochemical networks.
\newblock {\em PLoS Comput Biol}, 5(1):e1000256, 01 2009.

\bibitem{Multiscale}
Hannes Hettling and Johannes~HGM van Beek.
\newblock Analyzing the functional properties of the creatine kinase system
  with multiscale {`}sloppy{'} modeling.
\newblock {\em PLoS Comput Biol}, 7(8):e1002130, 08 2011.

\bibitem{VDP}
Balthasar van~der Pol.
\newblock On relaxation-oscillations.
\newblock {\em The London, Edinburgh and Dublin Phil. Mag. \& J. of Sci.},
  2(7):978--992, 1927.

\bibitem{Jones}
Christopher~K.R.T. Jones and Alexander I.~Khibnik (Eds.).
\newblock {\em Multiple-time-scale dynamical systems}.
\newblock Springer, 2000.

\bibitem{GrasmanBook}
Johan Grasman.
\newblock {\em Asymptotic Methods for Relaxation Oscillations and
  Applications}.
\newblock Springer Press, 1987.

\bibitem{Strogatz}
Steven~H. Strogatz.
\newblock {\em Nonlinear Dynamics and Chaos: With Applications to Physics,
  Biology, Chemistry and Engineering}.
\newblock Westview Press, 2001.

\bibitem{Note2}
Incidentally, the set of equations (\ref {vdpeqns}) can also be considered a
  special case of the FitzHugh-Nagumo model (\protect \textit {See} R.
  FitzHugh. Impulses and Physiological States in theoretical models of nerve
  propagation \protect \textit {Biophys J.}, 1(6):445, 1961) introduced three
  decades later as simplification of the Hodgkin-Huxley equations of neuronal
  spikes in the squid giant axons, and is sometimes referred to as the
  Bonhoeffer-van der Pol model.

\bibitem{Katok}
Anatole Katok and Boris Hasselblatt.
\newblock {\em Introduction to the Modern Theory of Dynamical Systems}.
\newblock Cambridge University Press, 1997.

\bibitem{Gucken}
John~M. Guckenheimer and Phillip Holmes.
\newblock {\em Nonlinear Oscillations, Dynamical Systems and Bifurcations of
  Vector Fields}.
\newblock Springer Press, 1983.

\bibitem{Unfold}
J.~Murdock.
\newblock {\em Normal Forms and Unfoldings for Local Dynamical Systems}.
\newblock Springer, New York, 2003.

\bibitem{Note3}
Perturbations distort the dynamics so that the attractor and its period change.
  We addressed these issues by setting the periods to unity, and by moving the
  initial conditions to the new attractor to remove any transients.
  Alternatively, if we fit data over many periods without making the said
  changes, the parameter combinations determining the period and phase would
  become stiff modes in our dynamics.

\bibitem{SloppyCell1}
Ryan~N. Gutenkunst, Jordan~C. Atlas, Fergal~P. Casey, Robert~S. Kuczenski,
  Joshua~J. Waterfall, Christopher~R. Myers, and James~P. Sethna.
\newblock Sloppycell \texttt{http://sloppycell.sourceforge.net}.
\newblock 2007.

\bibitem{SloppyCell2}
Christopher~R. Myers, Ryan~N. Gutenkunst, and James~P. Sethna.
\newblock Python unleashed on systems biology.
\newblock {\em Computing in Science Engineering}, 9(3):34 --37, 2007.

\bibitem{Note5}
We understand this as sloppiness as arising due to the generalized
  interpolation argument~\cite {MarkLongPaper}.

\bibitem{NR}
William~H. Press, Saul~A. Teukolsky, William~T. Vetterling, and Brian~P.
  Flannery.
\newblock {\em Numerical Recipes 3rd Edition: The Art of Scientific Computing}.
\newblock Cambridge University Press, New York, NY, USA, 3 edition, 2007.

\bibitem{Lipson}
Josh Bongard and Hod Lipson.
\newblock Automated reverse engineering of nonlinear dynamical systems.
\newblock {\em Proceedings of the National Academy of Sciences},
  104(24):9943--9948, 2007.

\bibitem{Comment}
Ricky Chachra, Mark~K. Transtrum, and James~P. Sethna.
\newblock Comment on {``}{S}loppy models{,} parameter uncertainty{,} and the
  role of experimental design{''}.
\newblock {\em Mol. BioSyst.}, 7:2522--2522, 2011.

\end{thebibliography}

\end{document}